\documentclass[11pt,fleqn]{article}

\usepackage{amssymb,latexsym, amsmath, epsfig, psfrag}

\newcommand{\halmos}{\hfill $\Box$}

\newcommand{\pr}[1]{{\rm {\mathbb P}}\{#1\}}
\newcommand{\expect}[1]{{\rm {\mathbb E}}\{#1\}}

\newcommand{\expecteps}[1]{{\rm {\mathbb E}}_{\nu_\varepsilon}\{#1\}}
\newcommand{\expectnu}[1]{{\rm {\mathbb E}}_\nu\{#1\}}
\newcommand{\preps}[1]{{\rm {\mathbb P}}_{\nu_\varepsilon}\{#1\}}
\newcommand{\prnu}[1]{{\rm {\mathbb P}}_\nu\{#1\}}

\newcommand{\beeq}{\begin{equation}}
\newcommand{\eneq}{\end{equation}}
\newcommand{\bear}{\begin{eqnarray}}
\newcommand{\enar}{\end{eqnarray}}
\newcommand{\bearno}{\begin{eqnarray*}}
\newcommand{\enarno}{\end{eqnarray*}}
\newtheorem{Th}{Theorem}[section]

\newtheorem{Cor}[Th]{Corollary}
\newtheorem{Lemma}[Th]{Lemma}
\newtheorem{Prop}[Th]{Proposition}

\usepackage{fullpage}

\begin{document}

\title{A large-deviations analysis of the $GI/GI/1$ SRPT queue}

\author{
Misja Nuyens$^{\ast}$ and  Bert Zwart$^{\dagger,\ddagger}$\\
 \\ $^{\ast}${\normalsize Department of Mathematics}\\
{\normalsize Vrije Universiteit Amsterdam}\\
{\normalsize De Boelelaan 1081, 1081 HV Amsterdam, The Netherlands}\\
{\normalsize{\tt mnuyens@few.vu.nl}, phone  +31 20 5987834, fax +31 20 5987653}
\vspace{5mm}\\
$^{\dagger}${\normalsize CWI}\\
{\normalsize P.O. Box 94079, 1090 GB Amsterdam, The Netherlands}
\vspace{5mm}\\
$^{\ddagger}${\normalsize Department of Mathematics \& Computer Science}\\
{\normalsize  Eindhoven University of Technology}\\
{\normalsize P.O. Box 513, 5600 MB Eindhoven, The Netherlands}\\
{\normalsize {\tt zwart@win.tue.nl}, phone +31 40 2472813,  fax +31 40 2465995}
\vspace{5mm}}

\date{\today}

\maketitle

\begin{abstract}

\noindent
We consider a $GI/GI/1$ queue with the shortest remaining processing time
discipline (SRPT) and light-tailed service times.
Our interest is focused on the tail behavior of the sojourn-time
distribution.
We obtain a general expression for its large-deviations decay rate.
The  value of this decay rate critically depends on 
whether there is mass in the endpoint of the service-time 
distribution or not.  An auxiliary priority queue, for which we obtain some new results, plays an
important role in our analysis. We apply our SRPT-results  
to compare  SRPT with FIFO from a
large-deviations point of view.
{\it Keywords \& Phrases:}  busy period, large deviations, priority queue,
shortest remaining processing time, sojourn time.\\
{\it Short title:} Large deviations for SRPT

\end{abstract}

\newpage

\section{Introduction}
\setcounter{equation}{0}

In queueing theory the shortest remaining processing time (SRPT)
discipline is famous, since it is known to minimize the mean queue
length and sojourn time over all work-conserving disciplines, see for example Schrage \cite{Schrage} and
Baccelli \& Br\'emaud \cite{BB}. Recent developments in communication
networks have led to a renewed interest
in queueing models with SRPT. For example,
Harchol-Balter {\em et al.}~\cite{HSBA03}
propose the usage of SRPT in web servers.
An important issue in such applications is the performance of SRPT for customers with a
given service time.
Bansal \& Harchol-Balter \cite{BH01} give some evidence
against the opinion that SRPT does not work well for large jobs.
They base their arguments on mean-value analysis.
Some interesting results on the mean sojourn time in heavy traffic
were recently obtained by Bansal \cite{B04} and Bansal \& Gamarnik
\cite{BG05}, who show that SRPT significantly outperforms FIFO if the
system is in heavy traffic. 

In the present paper we approach SRPT from a large-deviations
point of view. We  investigate the probability of a long sojourn time,
assuming that service times are light-tailed.
For heavy-tailed (more precisely, regularly varying)
service-time distributions,
N\'u\~nez-Queija \cite{Nunez} has shown that the tail of the sojourn-time
 distribution $\pr{V_{SRPT}>x}$ and the tail  of the service-time
distribution $\pr{B>x}$ coincide up to a constant. This appealing
property is shared by several
other preemptive service disciplines, for example by Last-In-First-Out (LIFO),
Foreground-Background (FB)
and Processor Sharing (PS);
see \cite{BBNZ03} for a survey. Non-preemptive service disciplines, like
FIFO, are known to behave worse: the tail of the sojourn time behaves like
 $ x\pr{B>x}$. This is the worst possible case, since it coincides with the
tail behavior of
a residual busy period; for details see again~\cite{BBNZ03}.

For light-tailed service times the situation is reversed. In a
fundamental paper, Ramanan \& Stolyar \cite{RS01} showed
that FIFO maximizes the decay rate (see Section \ref{s2} for a precise definition)
of the
sojourn-time distribution over all work-conserving
service disciplines. Thus, from a large-deviations point of view,
FIFO is optimal for light-tailed service-time distributions. Since for
any work-conserving service discipline the sojourn time is bounded by
a residual busy period, the decay rate of the residual busy period is again
the worst possible.
Recently, it has been shown that this worst-case decay-rate behavior of
the sojourn time  is  exhibited  under LIFO,
FB \cite{MN05}, and, under an additional assumption, PS \cite{MZ05}.

The present paper shows that a similar result holds for both non-preemptive
and preemptive SRPT,
under the assumption that the service-time distribution has no
mass at its right endpoint. Thus, for many light-tailed service-time
distributions, as for example phase-type service times, large sojourn times
are much more likely under SRPT than under FIFO. The derivation
of this result is based upon  a simple probabilistic argument; see Section \ref{4.1}.

The case where there is mass at the right endpoint of the service-time
distribution
 may be considered to be a curiosity; however, from a theoretical point
of view, it  actually turns out to be the
most interesting case. The associated analysis,
carried out in Section \ref{4.2}, is based on a  relation
with a $GI/GI/1$ priority queue. Since we could not find large-deviations
results in the literature (an in-depth treatment of the $M/G/1$ priority
queue is provided by Abate \& Whitt \cite{AW97}),
we  analyze this $GI/GI/1$ priority queue in Section \ref{s3}.
Another noteworthy feature of this case is that the resulting decay rate
is strictly larger than the one under LIFO, but strictly smaller
than under FIFO (with the exception of  deterministic  service times,
for which the FIFO decay rate is attained).
A similar result was recently shown in Egorova {\em et al.}~\cite{EZB05}
for the $M/D/1$ PS queue.
However, in general examples of service disciplines that exhibit this 
``in-between'' behavior are rare; see Section 5.1 of this paper for an overview.

Our results on SRPT suggest that, from a large-deviations point of view, it is not advisable
to switch from FIFO to SRPT.
However, in Section 6 we show that this suggestion should be handled with
care. Specifically, we investigate the decay rate of the {\em conditional} sojourn time, i.e.,
the sojourn time of a customer with service time $y$. We show that there exists
a critical service time $y^*$ such that SRPT is better than FIFO for service times below $y^*$ and
worse for service times larger than $y^*$. A performance indicator is the fraction of customers 
with service time exceeding $y^*$. We show that this fraction is close to zero
for both low and high loads;
numerical experiments suggest that this fraction is still very small
for moderate values of the load.

This paper is organized as follows. Section \ref{s2} introduces notation and
states some preliminary results. In particular, the decay rates of the
workload and busy period are derived in complete generality. Section 3
treats a
two-class priority queue with renewal input and investigates the
tail behavior of the low-priority waiting time.
The results on SRPT are presented in Section \ref{s4}.
Section \ref{s5} treats various implications of the results in Sections \ref{s3} and \ref{s4}.
First, we  compare our results with the decay rates for LIFO and FIFO, and
show that the decay rate of the sojourn time under SRPT is strictly in between these two
if the service-time distribution has mass at its right endpoint. We then treat the special case of Poisson arrivals; in particular
we show that our results for the priority queue agree with those of Abate \& Whitt \cite{AW97}.
In addition, we consider the behavior of the decay rates in heavy traffic.
Conditional sojourn times are investigated in Section \ref{s6}. We summarize our results and propose  directions for further
research in Section \ref{s7}.

\section{Preliminaries: workload and busy period}\label{s2}
\setcounter{equation}{0}

In this section we introduce the notation and derive two
 preliminary results.
We consider a stationary, work-conserving
$GI/GI/1$ queue, with the server
working at unit speed.
Generic inter-arrival and service times are denoted by $A$ and $B$.
To avoid trivialities, we assume that $\pr{B>A}>0$ (otherwise there would be
no delays).
Define the system load $\rho=\expect{B}/\expect{A}<1$.
Since $\rho<1$, the workload process is positive recurrent and the busy
period $P$ has finite mean.
The moment generating function of a random variable $X$ is  denoted by
 $\Phi_X(s)=\expect{{\rm e}^{sX}}$.
Throughout the paper we assume that $B$ is light-tailed, i.e., that $\Phi_B(s)$ is finite in a
neighborhood of $0$.
Let $W$ be the workload seen by a customer upon arrival in steady state.
This workload
 coincides with the FIFO waiting time. Furthermore, let $W^y$ be the steady-state workload
on arrival epochs in the $GI/GI/1$ queue with service times $B^y=BI(B<y)$.
Let
$P^y$ denote the busy period in such a queue.
Our first preliminary result concerns the logarithmic tail asymptotics for $W$.

\begin{Prop}
\label{logW}
As $x\rightarrow\infty$, we have that $\log\pr{W>x} \sim - \gamma_{w} x$,
with
\beeq \label{eerste}\gamma_w=\sup\{ s: \Phi_A(-s)\Phi_B(s)\leq 1 \}.\eneq
\end{Prop}

We call $\gamma_w$ the {\em decay rate} of $W$. Generally, for any random
variable $U$, we call $\gamma_u$ the decay rate of $U$ if for $x\rightarrow\infty$,
\[\log \pr{U>x} =-\gamma_ux+{\rm o}(x).\]

If $\Phi_A(-\gamma_w)\Phi_B(\gamma_w)=1$,
several proofs of Proposition \ref{logW} are available, see e.g.~Asmussen \cite{Asmussen03}, Ganesh {\em et al.}~\cite{GOW} and Glynn \& Whitt \cite{GW94}.
 We believe that the result in its
present generality is known as well, but could not find a reference.
For completeness,  a short proof is included here. \\ \\
{\bf Proof of Proposition \ref{logW}} \\
The upper bound follows from a famous result of Kingman \cite{Kingman}:
\[
\log \pr{W>x} \leq -\gamma_w x.
\]
For the lower bound we use a truncation argument.
 From Theorem XIII.5.3 of
\cite{Asmussen03}
 (the condition of that theorem is easily seen to be satisfied for bounded service times), it follows that
\[
\log \pr{W^y>x} \sim -\gamma_w^yx,
\]
with $\gamma_w^y=\sup\{s: \Phi_A(-s)\Phi_{B^y}(s)\leq 1\}$.
Consequently, since $\pr{W>x} \geq \pr{W^y>x}$,
\[
\liminf_{x\rightarrow\infty} \frac{1}{x}\log\pr{W>x} \geq -\gamma_w^y.
\]
Since $\Phi_{B^y}(s)$ is increasing in $y$, and $\Phi_{B^y}(s)$ converges to $\Phi_{B}(s)$ as $y\to\infty$, the decay rate  $\gamma_w^y$ is decreasing
in $y$, and converges to a limit $\gamma_w^*\geq \gamma_w$.
Since $\gamma_w^*\in [0,\gamma_w^y]$ for any $y$, and 
$\Phi_A(-s)\Phi_{B^y}(s)$ is convex in $s$ and 
has a negative derivative in 0, we have
$\Phi_A(-\gamma_w^*)\Phi_{B^y}(\gamma_w^*)\leq 1$ for all $y$.
Consequently,
\[
\Phi_A(-\gamma_w^*)\Phi_{B}(\gamma_w^*)=
\lim_{y\rightarrow\infty} \Phi_A(-\gamma_w^*)\Phi_{B^y}(\gamma_w^*)\leq 1.
\]
This implies that $\gamma_w^*\leq \gamma_w$, so that $\lim_{y\to\infty} \gamma_w^y=\gamma_w^*=\gamma_w$. This yields the desired
lower limit.
\halmos \\

We continue by deriving
an expression for the decay rate $\gamma_p$ of the busy period $P$.
Sufficient conditions for precise asymptotics of $\pr{P>x}$, which are of the
form
$Cx^{-3/2}{\rm e}^{-\gamma_px}$, are given in Palmowski \& Rolski \cite{PR04}.
These asymptotics follow
from a detailed analysis, involving a change-of-measure argument.
We show that logarithmic asymptotics (which are of course implied by  precise asymptotics) can be given without any further assumptions.

\begin{Prop}\label{busyperiodprop}
As $x\rightarrow\infty$, we have $\log \pr{P>x} \sim - \gamma_p x,$
with
\beeq \label{gammap}\gamma_p=\sup_{s\geq 0} \{s-\Psi(s)\},\eneq
and $\Psi(s)=-\Phi_A^{-1} \left(\frac 1{\Phi_B(s)}\right)$.
\end{Prop}
{\bf Proof}\\
We first derive an upper bound.
 Let $X(t)$ be the amount of work offered to the
queue in the interval $[0,t]$.
In Lemma 2.1 of Mandjes \& Zwart \cite{MZ05}
it is shown that for each $s\geq 0$,
\beeq
\label{psi}
\Psi(s)=\lim_{t\rightarrow\infty} \frac 1t \log \expect{{\rm e}^{sX(t)}}.
\eneq
Using the Chernoff bound, we have for all $s\geq 0$,
\[
\pr{P>t} \leq \pr{X(t)>t} \leq {\rm e}^{-st+ \log \expect{\exp\{sX(t)\}}}.
\]
Consequently,
\[
\limsup_{t\rightarrow\infty} \frac 1t \log \pr{P>t} \leq -s +
\limsup_{t\rightarrow\infty} \frac {1}{t}\log \expect{\exp\{sX(t)\}}
=-(s-\Psi(s)).
\]
Minimizing over $s$  yields the upper bound for $\pr{P>t}$.
We now turn to the lower bound, for which we again use a truncation argument.
First, note that
\[
\pr{P>x} \geq \pr{P^y>x}.
\]
For truncated service times, the assumptions in \cite{PR04}
for the exact asymptotics (cf.\ Equation (33) in \cite{PR04}) are  satisfied, and we have,
with obvious notation,
\[
\liminf_{x\to\infty} \frac{1}{x}\log\pr{P>x}\geq
\lim_{x\rightarrow\infty} \frac{1}{x}\log\pr{ P^y>x}=
-\sup_{s\geq 0}\{s-\Psi^y(s)\}=-\gamma_p^y.
\]
So to prove the theorem, it suffices to  show that $\gamma_p^y\to \gamma_p$ for $y\to\infty$.
Define $f^y(s)=s-\Psi^y(s)$. It is clear that
$f^y(s)\rightarrow f(s)= s-\Psi(s)$ pointwise as $y\rightarrow\infty$ and that
$f^y(s)$ is decreasing in $y$.
Consequently, we have that the limit of $\gamma_p^y$ for $y\rightarrow\infty$ exists and
that
\[
\gamma_p^*=\lim_{y\to\infty} \gamma_p^y=\lim_{y\rightarrow\infty} \sup_{s\geq 0} f^y(s)
\geq \sup_{s\geq 0} f(s)=\gamma_p.
\]
It remains to show that the reverse inequality holds. For this, we use an argument similar to one
in the proof of Cram\'ers theorem (cf.~Dembo \& Zeitouni \cite{DZ98}, p.\ 33).
Take $y_0$ such that $\pr{B^y>A}>0$ for $y>y_0$.
Then there exist $\delta,\eta>0$ such that $\pr{B^y-A\geq \delta}
\geq \eta>0$ for $y\geq y_0$.
Hence, for $y\geq y_0$,
\[ \Phi_{B^y}(s)\Phi_A(-s)= \expect{{\rm e}^{sB^y}}
\expect{{\rm e}^{-sA}}=\expect{{\rm e}^{s(B^y-A)}} \geq \eta e^{s\delta}.\]
For $s$ large enough, we now have
\[ \Phi_A(-s)\geq \frac{1}{\Phi_{B^y}(s)}. \]
Since $\Phi_A^{-1}(s)$ is increasing in $s$, we find that for $s$ and $y$ large enough,
\[ s+\Phi_A^{-1}\Big(\frac{1}{\Phi_{B^y}(s)}\Big)\leq s+\Phi_A^{-1}\big( \Phi_A(-s)\big)=0.\]
Since $\Phi_A^{-1}\left(1/\Phi_{B^y}(s)\right)$ is decreasing in $y$
and is continuous in $s$, we see that for $y\geq y_0$
the level sets $L_y=\{s: f^y(s)\geq \gamma_p^*\}$
are compact. Moreover, since $f^y(s)$ is decreasing in $y$,
the level sets are nested with respect to $y$. Consequently, the intersection of the level
sets $L_y$ contains at least one element, say $s_0$.
By the definition of $s_0$, we have $f^y(s_0)\geq \gamma_p^*$ for every $y$. Thus, since $f^y$ converges pointwise,
\[
\gamma_p=\sup_{s\geq 0} f(s) \geq f(s_0)= \lim_{y\rightarrow\infty} f^y(s_0) \geq \gamma_p^*.
\]
We conclude that $\gamma_p^y \rightarrow \gamma_p$ as $y\rightarrow\infty$, which completes the proof.\halmos

\section{The $GI/GI/1$ priority queue}\label{s3}
\setcounter{equation}{0}

In this section, we consider the following $GI/GI/1$ two-class priority queue.
Customers arrive according to a renewal process with generic inter-arrival time
$A$. An arriving customer is of class 1 with probability $p$,
in which case he has service time $B_1$. Customers of class 2 have service time
 $B_2$.
Class-1 customers have priority over class-2 customers. We assume that $0<p<1$,
and that $p\expect{B_1}+(1-p)\expect{B_2}<\expect{A}$,
which ensures that the priority queue is stable.
We are interested in the steady-state waiting time $W_2$ of a class-2 customer, that is,
the time a class-2 customer has to wait before he enters service for the first
 time. Note that $W_2$  is independent of whether the priority mechanism is preemptive or
 not.

Let $N_1(t)$ be the renewal process generated by the arrivals of the class-1 customers, i.e.,
$N_1(t)=\max \{ n: A_{1,1}+\cdots + A_{1,n} \leq t\}$. Here $A_{1,i}$ is the time
between the arrival of the $(i-1)$-st and $i$-th customer. A generic class-1
inter-arrival time is denoted by $A_1$.
Note that $A_1$ is a geometric sum of ``original'' inter-arrival times $A$:
\[
\Phi_{A_1}(s)=\expect{{\rm e}^{sA_1}}= \sum_{n=0}^\infty p(1-p)^n
\Phi_A(s)^{n+1}= \frac{p\Phi_A(s)}{1-(1-p)\Phi_A(s)}.
\]
Define
\[
X_1(t)= \sum_{i=1}^{N_1(t)} B_{1,i}.
\]
Hence, $X_1(t)$ is the amount of work of type $1$  that has arrived in the system by time $t$. Let $P_{1}$ be a generic busy period of class $1$ customers.
Finally, let $P_1(x)$ be a busy period of class-1 customers with an
{initial customer of size} $x$, so
\[
P_1(x)\stackrel{d}{=} \inf \{t\geq 0: x+X_1(t)\leq t \}.
\]
Denoting the total workload in the queue at arrivals again by $W$ (cf.\ Section 2),
we have the following fundamental identity:
\beeq \label{p1w}
W_2 \stackrel{d}{=} P_1(W),
\eneq
where $W$  and $\{P_1(x),x\geq 0\}$ are independent.
This identity holds since, using a discrete-time version of PASTA, $W$ is also
the workload as seen by an arriving customer of class 2.
Set
\beeq \label{psi1} \Psi_1(s)=-\Phi_{A_1}^{-1}\left(\frac 1{\Phi_{B_1}(s)}\right).\eneq
 The main result of this section is the following.

\begin{Th}\label{w2thm}
As $x\rightarrow\infty$, we have $\log \pr{W_2>x} \sim - \gamma_{w_2} x,$
with
\beeq\label{drw2}
\gamma_{w_2}=\sup_{s\in [0,\gamma_w]} \{s-\Psi_1(s)\}.
\eneq
\end{Th}

Before we give a proof of this theorem, we first describe some heuristics,
starting from $W_2\stackrel{d}{=}P_1(W)$.
The most likely way for $W_2$ to become large (i.e., $W_2>x$) involves 
a combination of two events: (i) $W$ is of the order $ax$ for some
constant $a\geq 0$; (ii) $P_1(ax)$ is of the order $x$.
Clearly, there is a trade-off: as $a$ becomes larger, scenario (i) become less likely,
while scenario (ii) becomes more likely. Thus, we need to find the optimal value
of $a$. For this we need to know the large-deviations decay rates associated with events (i) and (ii).
The decay rate of event (i) is simply
$a\gamma_w$. To obtain the decay rate of event (ii), note that
\[
\pr{P_1(ax)>x} \approx \pr{X_1(x)>(1-a)x}.
\]
One can show that the RHS  probability has decay rate $\sup_{s\geq 0}\{(1-a)s-\Psi_1(s)\}$.
Thus, the optimal value of $a$, and the decay rate $\gamma_{w_2}$,
can be found by optimizing the expression
\[
\inf_{a\geq 0}\{a\gamma_w+ \sup_{s\geq 0}[(1-a)s+\Psi_1(s)]\}.
\]
It is possible to show that the value of this program coincides with
$\sup_{s\in [0,\gamma_w]} \{s-\Psi_1(s)\}$. Moreover, the optimal value of $a$ is $0$ if the optimizing argument
of $s-\Psi_1(s)$ is strictly less than $\gamma_w$, and  it is $1-\Psi_1'(\gamma_w)$ if
$\sup_{s\in [0,\gamma_w]} \{s-\Psi_1(s)\}=\gamma_w-\Psi_1(\gamma_w)$.
In the proof below, we only use these heuristics to ``guess'' the correct value of $a$.

Note that the two cases $a>0$ and $a=0$ correspond to two qualitatively different scenarios leading to a large value of $W_2$. If $a=0$, then the customer sees a ``normal" amount of work upon arrival, while $a>0$ results in a workload of the order $ax$ at time $0$. This distinction between two different scenarios is typical in priority queueing, see
 Abate \& Whitt \cite{AW97} and Mandjes \& Van Uitert \cite{MU05} for more discussion.
\\

\noindent
{\bf Proof} \\
We start with the upper bound. Using the Chernoff bound, we find that for $s\geq 0$,
\[
\pr{W_2>x}=\pr{P_1(W)>x} \leq \pr{W+X_1(x)-x>0} \leq \expect{ {\rm e}^{sW} }
{\rm e}^{-xs}\expect{{\rm e}^{sX_1(x)}}.
\]
Using (\ref{psi}) with $X(t)$ replaced by $X_1(t)$, we see that for all $s\in[0, \gamma_w)$,
\[
\limsup_{x\rightarrow\infty} \frac 1x \log \pr{W_2>x} \leq -[s-\Psi_1(s)].
\]
The proof of the upper bound is completed by minimizing over $s$, and noting that
$\sup_{s\in [0,\gamma_w)} \{s-\Psi_1(s)\}=\sup_{s\in [0,\gamma_w]} \{s-\Psi_1(s)\}$.

We now turn to the lower bound.
 From the proof of Proposition \ref{busyperiodprop}, we see that
$P_1$ has decay rate $\gamma_{p_1}=\sup_{s\geq 0} \{s-\Psi_1(s) \}$.
Let $s_1$ be the unique optimizing argument.
In addition,
let $r$ be the probability that  at the arrival of a class-2 customer to the steady state queue
at  least one customer of type 1 is waiting.
It is obvious that $r>0$. Since $P_1(W)\geq_{st} P_1$ on this event, we see that
\[
\pr{P_1(W)>x} \geq r\pr{P_1>x},
\]
which by (\ref{p1w}) implies that
\[
\liminf_{x\rightarrow\infty} \frac 1x \log \pr{W_2>x}\geq -\gamma_{p_1}.
\]
Thus, if $s_1\leq \gamma_w$, we can conclude from this and the upper bound that
\[
\lim_{x\rightarrow\infty} \frac 1x \log \pr{W_2>x}= -\gamma_{p_1}.
\]
What remains is to consider the case $s_1>\gamma_w$.
Since the concave function $s-\Psi_1(s)$ is increasing between $0$ and $s_1$,
we see that
$\sup_{s\in [0,\gamma_w]} \{s-\Psi_1(s)\}=\gamma_w-\Psi_1(\gamma_w)$.
Thus, to complete the proof of the theorem, it suffices to show that
\beeq
\label{toprove}
\liminf_{x\rightarrow\infty} \frac 1x \log \pr{W_2>x}\geq -[\gamma_w-\Psi_1(\gamma_w)].
\eneq
Note that for any $a>0$,
\beeq \label{w2:1}
 \pr{W_2>x} \geq  \pr{W>ax} \pr{P_1(ax)>x}
= {\rm e}^{-a \gamma_w x +{\rm o}(x)} \pr{P_1(ax)>x}.
\eneq
Combining (\ref{w2:1}) and Lemma \ref{changeofm} below, we see that
by taking $a=1-\Psi_1'(\gamma_w)$,
\[
\log \pr{W_2>x} \geq  -a \gamma_w x +{\rm o}(x)+\log  \pr{P_1(ax)>x}
= -x(\gamma_w-\Psi_1(\gamma_w))+{\rm o}(x),
\]
which coincides with (\ref{toprove}), as was required.
\halmos\\

We now provide the result that was quoted in the proof above.
\begin{Lemma}\label{changeofm}
 Set $a=1-\Psi_1'(\gamma_w)$. If $\gamma_w<s_1$, then
\[
\log \pr{P_1(ax)>x}\geq  -x(\gamma_w(1-a)-\Psi_1(\gamma_w))+{\rm o}(x).
\]
\end{Lemma}

\noindent
{\bf Proof}\\
To prove the lemma we use a change-of-measure argument.
Define a probability measure $\prnu{\cdot}$ for $\nu\geq 0$ such
that
\bearno
\prnu{A_{1,i} \in {\rm d}x} &=& {\rm e}^{-\Psi_1(\nu)x} \pr{A_{1,i} \in {\rm d}x}/\Phi_{A_1}(-\Psi_1(\nu)),
\hspace{1cm} i\geq 1,\\
\prnu{B_{1,i} \in {\rm d}x} &=& {\rm e}^{\nu x} \pr{B_{1,i} \in {\rm d}x}/\Phi_{B_1}(\nu), \hspace{1cm} i\geq 1.
\enarno
Choose $\nu=\nu_\varepsilon$ such that
\[
\Psi_1'(\nu_\varepsilon)=\frac{\expectnu{B_i}}{\expectnu{X_i}}=
\frac{\Phi_{B_1}'(\nu_\varepsilon)}{\Phi_{B_1}(\nu_\varepsilon)}\left/
\frac{\Phi'_{A_1}(-\Psi(\nu_\varepsilon))}{\Phi_{A_1}(-\Psi(\nu_\varepsilon))}\right.=
1-a+\varepsilon,\qquad  \varepsilon<a.
\]
We denote this probability measure by $\preps{\cdot}$.
The drift under this new measure is $1-a+\epsilon$, making the event 
$\{P_1(ax)>x\}$ extremely likely
for large $x$.
Note that $\nu_0=\gamma_w$, by the definition of $a$, and since $\Psi_1'(s)$ is strictly increasing. \\

Let ${\cal F}_n$ be the Borel $\sigma$-algebra generated by
$A_{1,1},\ldots, A_{1,n}, B_{1,1},\ldots,B_{1,n}$. Define $S_n^{A_1}=A_{1,1}+\ldots + A_{1,n}$ and
$S_n^{B_1}=B_{1,1}+\ldots +B_{1,n}$.
Note that $\bar N_1(x):=N_1(x)+1$ is a stopping time w.r.t.~the filtration $({\cal F}_n).$ Furthermore, note
that the event $\{P_1(ax)>x\}$ is ${\cal F}_{\bar N(x)}$-measurable.
Finally, note that
for every $\varepsilon>0$ small enough, the process 
$1/M_n^\varepsilon, n\geq 1$, with
\[
M_n^\varepsilon=\exp\{\Psi_1(\nu_\varepsilon)S_n^{A_1}-\nu_\varepsilon S_n^{B_1}\},
\]
is a martingale w.r.t.~${\cal F}_n$ under $\pr{\cdot}$,
since the definition of $\Psi_1$ ensures that $\Phi_{A_1}(-\Psi_1(\nu_\varepsilon))\Phi_{B_1}(\nu_\varepsilon)=1$.
Thus, we have the following fundamental identity (see for example~Theorem XIII.3.2 in \cite{Asmussen03}):
\[
\pr{P_1(ax)>x} = \expecteps{M_{{\bar N}_1(x)}^\varepsilon I(P_1(ax)>x)}.
\]
Furthermore, we have for any event ${\cal S} \subseteq {\cal F}_{\bar N_1(x)}$,
\beeq
\label{lb1}
\pr{P_1(ax)>x} \geq \expecteps{M_{{\bar N}_1(x)}^\varepsilon I(P_1(ax)>x) I({\cal S})}.
\eneq
Take here
\[
{\cal S}\equiv {\cal S}_{\varepsilon}:=\left\{ S^{B_1}_{N_1(x)}\leq (1-a+\varepsilon)x \right\}.
\]
Note that  $S_{N_1(x)+1}^{A_1}> x$ by definition and apply
the definition of ${\cal S}_\varepsilon$ to obtain from (\ref{lb1})  the following
lower bound for $\pr{P_1(ax)>x}$:
\[
\liminf_{x\to\infty}\frac 1x \log\pr{P_1(ax)>x}\geq
-\nu_\varepsilon(1-a+\varepsilon)+\Psi_1(\nu_\varepsilon)+
\liminf_{x\to\infty}\frac 1x \log \preps{P_1(ax)>x,{\cal S}_\varepsilon}.
\]
By the law of large numbers, we have that
$\preps{P_1(ax)>x; {\cal S}_\varepsilon}$ is bounded away from zero, uniformly in $x$ for every $\varepsilon>0$.
Consequently,
\[
\liminf_{x\to\infty}\frac 1x \log\pr{P_1(ax)>x}\geq
-\nu_\varepsilon(1-a+\varepsilon)+\Psi_1(\nu_\varepsilon).
\]
Now let $\varepsilon\downarrow 0$. This yields
\[
\liminf_{x\to\infty}\frac 1x \log\pr{P_1(ax)>x} \geq  -\gamma_w(1-a)+\Psi_1(\gamma_w),
\]
and the statement of the lemma follows. \halmos \\

 From Theorem \ref{w2thm} we can deduce the decay rate of the sojourn time $V_2$ of class-2 customers.
This turns out to be the same for both preemptive and non-preemptive service.

\begin{Th} As $x\to\infty$, we have $\log \pr{V_2>x}\sim-\gamma_{w_2}x$, where $\gamma_{w_2}$ is as in (\ref{drw2}).
\end{Th}
{\bf Proof}\\
For the non-preemptive case, we have $V_2=W_2+B_2$, where $W_2$ and $B_2$ are independent.
Since the decay rate of $B_2$ is larger than $\gamma_{w_2}$, and since the decay rate of a sum
of independent random variables is equal to the smallest decay rate (see for example\ \cite{MN05} for a
short proof), the result for this case follows immediately.
In the preemptive case, we use that
\[
V_2\stackrel{d}{=} P_1(W+B_2)\geq_{st} W_2,
\]
which gives us the lower bound.
The upper bound follows the same lines of the proof of Theorem \ref{w2thm} and
noting that $\expect{{\rm e}^{\gamma_w B_2}}<\infty$.
\halmos

\section{Shortest Remaining Processing Time}\label{s4}
\setcounter{equation}{0}
In this section we present our results on the sojourn time under the SRPT discipline.
Define  $V_{SRPT}$ as the steady-state sojourn time of a customer under
the preemptive SRPT discipline.
Further, define the right endpoint $x_B$ by  $x_B=\sup\{ x: \pr{B> x}>0\}$.
When it comes to determining the decay rate of $V_{SRPT}$, it turns out to be crucial whether
\beeq
\label{xb}
\pr{B=x_B}=0,
\eneq
or not.
In the first subsection, we show that if (\ref{xb}) holds, then the decay rate of $V_{SRPT}$ is
equal to $\gamma_p$, the decay rate of the busy period $P$.
 If (\ref{xb}) does not
hold, the situation is more complicated. In that case we  use
  the results of the previous section to show that the decay rate of $V_{SRPT}$ is equal to $\gamma_{w_2}$,
where $W_2$ is the waiting time in a certain auxiliary priority queue.
This is the subject of the second subsection.
We also show that  for the non-preemptive SRPT discipline the same results hold.

\subsection{No mass at the right endpoint}\label{4.1}

In this section we prove the following theorem.

\begin{Th}\label{nomassinendpoint}
Suppose that $\pr {B=x_B}=0$.
Then  $\log \pr{V_{SRPT}>x} \sim  -\gamma_p x$ for $x\to\infty$,  with $\gamma_p$  as in (\ref{gammap}).
\end{Th}
{\bf Proof}\\
Let $V_{SRPT}$ be the sojourn time of a tagged customer with service time $B$.
Since $V_{SRPT}\leq P^*$, where $P^*$ is the residual busy period $P(W)$,
and since for light tails the decay rate of $P^*$ coincides with that of $P$ (this follows from Lemma 3.2 in \cite{AW97}),
we see that
\[
\limsup_{x\rightarrow \infty} \frac 1x \log \pr{V_{SRPT}>x} \leq -\gamma_p.
\]
Thus, it suffices to show that the corresponding result holds for the lower limit.
 For this, we construct a lower bound for $\pr{V_{SRPT}>x}$.
Assume first that $x_B=\infty$. Let $A$ be the last inter-arrival time before the tagged customer arrives,
 $B_0$ be the service time of that customer, and  $a$ be such that $\pr{A<a}>0$.
Then, for all $y$,
\bearno
\pr{V_{SRPT}>x} &\geq& \pr{V_{SRPT}>x; B> y, A<a, B_0\leq y} \\
                &\geq& \pr{A<a} \pr{B>y} \pr{B_0\leq y} \pr{P^{y-a}>x}.
\enarno
The last inequality holds since conditional on $A<a, B>y$ and $B_0\leq y$, the tagged customer
has to wait at least for the sub-busy period generated by the customer that
arrived before him,  and this sub-busy period is stochastically larger than  $P^{y-a}$.
Since $\pr{A<a}\pr{B>y}>0$, and $\pr{B_0\leq y}>0$ for $y$ large enough, we have that
\[
\liminf_{x\rightarrow \infty} \frac 1x \log \pr{V_{SRPT}>x} \geq -\gamma_p^{y-a}
\]
for $y$ large enough. Letting $y\rightarrow x_B=\infty$, we obtain $\gamma_p^{y-a}\rightarrow \gamma_p$,
as in the proof of Proposition 2.2.

If $x_B<\infty$, the above proof can be modified in a straightforward way if
$\pr{A<a}>0$ for all $a>0$.
However, this may not be the case in general and therefore we have to
make a more involved construction.
By definition of $x_B$, there exists a
decreasing sequence $(\varepsilon_n)$ such that
$\pr{x_B-\varepsilon_n< B< x_B-\varepsilon_n/2}>0$ for all $n$,
and $\varepsilon_n\to 0$ as $n\to\infty$.
Since $\pr{B>A}>0$, we can assume that $\varepsilon_1$ is such that
$\pr{A<x_B-2\varepsilon_1}>0$.
Let $R_n$ be the event that the last $\lfloor x_B/\varepsilon_n\rfloor$ customers
that arrived before the tagged customer had a service time in the interval
 $[x_B-\varepsilon_n,x_B-\varepsilon_n/2]$, and that
 the last $\lfloor x_B/\varepsilon_n\rfloor $ inter-arrival times were smaller than
 $x_B-2\varepsilon_n$. By definition of $\varepsilon_n$, we have $\pr{R_n}>0$ for all $n$.

Furthermore,
by the SRPT priority rule,  we see by induction that on the event $R_n$,
after the $k$th of the last $n$ inter-arrival times, there is a customer with remaining service time larger than $k\varepsilon_n$.
Hence, at the arrival of the tagged customer, there is a customer in the system with remaining service
time in the interval $[x_B-\varepsilon_n,x_B-\varepsilon_n/2]$.
If the tagged customer has service time $B>x_B-\varepsilon_n/2$,
his sojourn time satisfies
$V_{SRPT} \geq P^{x_B-\varepsilon_n}$ on $R_n$.
Consequently, for all $n\in\mathbb{N}$,
\[
\pr{V_{SRPT}>x} \geq \pr{R_n}\pr{B>x_B-\varepsilon_n/2} \pr{P^{x_B-\varepsilon_n}>x}.
\]
This implies
\[
\liminf_{x\rightarrow \infty} \frac 1x \log \pr{V_{SRPT}>x} \geq 
-\gamma_p^{x_B-\varepsilon_n}.
\]
Letting $n\to\infty$, and hence $\varepsilon_n\downarrow 0$, we get
$\gamma_p^{x_B-\varepsilon_n}\rightarrow \gamma_p$, as before. This completes the proof.~\halmos\\

The property that the decay rate of the sojourn time is equal to that of the busy period is shared by a number of disciplines, see Section \ref{comparison}. Further, we remark that for light tails, $\gamma_p$ is the smallest possible decay rate for the sojourn time in the class of all work-conserving disciplines: the sojourn time is bounded above by the residual busy period $P^*$, and for light-tailed service times $P^*$  has decay rate $\gamma_p$ (cf.\  Lemma 3.2 in \cite{AW97}).

\subsection{Mass at right endpoint}\label{4.2}

If there is mass at the right endpoint $x_B$ of the service-time distribution, then the tail behavior of  $V_{SRPT}$ is more complicated.
To obtain the decay rate of $V_{SRPT}$ for this case, we
 identify the SRPT queue with  the following two-class priority queue.
Let the customers of class 1 be the customers with service time strictly less
 than $x_B$. Then  $B_2=x_B$ and $B_1$ is such that
\beeq \label{b1} \pr{B_1\leq x}=\pr{B\leq x\mid B<x_B}, \qquad x\geq 0.\eneq

\begin{Th}\label{massinendpoint}
Suppose that $\pr{B=x_B}>0$. Then $\log \pr{V_{SRPT}>x} \sim  -\gamma_v x$ for $x\to\infty$,
with
\beeq
\label{drvsrpt}
\gamma_v=\sup_{s\in [0,\gamma_w]} \{s-\Psi_1(s)\},
\eneq
where $\Psi_1$ is as in (\ref{psi1}), and $B_1$ is as in (\ref{b1}).
\end{Th}
{\bf Proof}\\
First, note that if $q=\pr{B=x_B}=1$, we have a $G/D/1$ SRPT queue, which has the
 same dynamics as a FIFO queue. Indeed we obtain $\Psi_1\equiv 0$, implying
 $\gamma_v=\gamma_w$, cf.~Proposition 2.1.
Assume therefore that $0<q<1$, let $V_{SRPT}$ be the sojourn time of a tagged customer with service time $B$, and write
\[\pr{V_{SRPT}>x}=q\pr{V_{SRPT}>x\mid B=x_B}+(1-q)\pr{V_{SRPT}>x\mid B<x_B}.
\]
 From the nature of the SRPT discipline, or a simple coupling argument, it is obvious that
\[\pr{V_{SRPT}>x\mid B<x_B}\leq \pr{V_{SRPT}>x\mid B=x_B}.\]
 Therefore, it  suffices
to consider the tail behavior of $\bar V_{SRPT}$, where
\[\pr{\bar V_{SRPT}\leq x}=\pr{V_{SRPT}\leq x \mid B=x_B}.\]
First, we note that $\bar V_{SRPT}$ is bounded from below by the time it
 takes until our tagged customer receives
service for the first time. A crucial observation is that this period
coincides with the low priority waiting-time $W_2$ defined in Section 3.
Second, note that $\bar V_{SRPT}$ is upper bounded by the sojourn
time $V_2$ of a class-2 customer in the above priority queue. Hence, we have $W_2\leq_{st}\bar  V_{SRPT} \leq_{st} V_2$.

Further,  $V_2$ satisfies
$V_2\stackrel{d}{=} P_1(W+x_B)$ for preemptive service;
for non-preemptive service, we have $V_2\stackrel{d}{=} P_1(W)+x_B$.
Since the logarithmic asymptotics of $W+x_B$ coincide with those of $W$,
we can mimic the proof of Theorem \ref{w2thm} to see that in both cases
the decay rate of $V_2$ coincides with that of $W_2$.
Hence the decay rate of $\bar V_{SRPT}$ is given by (\ref{drw2}), and the proof is completed.
\halmos\\

The intuition of how $V_{SRPT}$ becomes large is the same as that of $W_2$ in Section 3.
In Section \ref{comparison} below, we show that if there is mass in the endpoint $x_B$, then the decay rate of the sojourn time lies strictly between the maximal value (obtained for FIFO) and the minimal value (LIFO).

\section{Complements}\label{s5}
\setcounter{equation}{0}

In the previous two sections we have derived expressions
for the decay rates $\gamma_{w_2}$ and $\gamma_v$.
In this section we derive some properties of these decay rates.
Specifically, in Section \ref{comparison} we compare $\gamma_{w_2}$ and $\gamma_{v}$ with $\gamma_w$
and $\gamma_p$. We show that for $q=\pr{B=x_B}\in(0,1)$, we always have
$\gamma_p<\gamma_{w_2}<\gamma_w$. Consequently, if $q\in (0,1)$,
we also find that $\gamma_v=\gamma_v(q)\in (\gamma_p,\gamma_w)$.
As explained in the introduction, this is a non-standard result.
We also indicate that $\gamma_v$ can take any value between $\gamma_p$
and $\gamma_w$, depending on the value of $q$.

Further, in Section \ref{possonian}  we specialize our expression of $\gamma_{w_2}$ to the case of Poisson arrivals. For priority queues, a quite involved expression for
the decay rate was given in \cite{AW97}. We show that this expression
can be simplified, and that it coincides with our expression of
$\gamma_{w_2}$.

Finally, in Section \ref{heavytraffic}, we derive heavy-traffic approximations
for $\gamma_{w_2}$.

\subsection{Comparison with other service disciplines} \label{comparison}

In this subsection, we compare the decay rates $\gamma_{w_2}$ and $\gamma_v$
with the decay rate of the sojourn time under FIFO and LIFO,
which respectively equal $\gamma_w$ and $\gamma_p$.

We first show
that for the priority queue described in Section \ref{s3}, 
the decay rate of $W_2$ is different from those of $P$ and $W$.

\begin{Prop}\label{diffdecayrates} Assume  $0<p<1$. Then $\gamma_p< \gamma_{w_2}<\gamma_w.$
\end{Prop}
{\bf Proof}\\
Since $\Psi_1(s) > 0$ for $s>0$,
we have by Theorem \ref{w2thm} that
\[\gamma_{w_2} = \sup_{s\in [0,\gamma_w]} \{s-\Psi_1(s)\}  <
\sup_{s\in [0,\gamma_w]} s  = \gamma_w.\]
To prove the inequality $\gamma_{w_2}>\gamma_p$, we provide a
different construction of the function $\Psi_1(s)$.
Let $B_p$ be a service time which is equal to $B_1$ with probability $p$ and
$0$ with probability $1-p$. It is clear that $\Phi_{B_p}(s)<\Phi_B(s)$.
The amount of work $X_1(t)$ generated by class-1 customers between
time $0$ and $t$ is the same in distribution as the amount of
work generated by the arrival process with inter-arrival times $A$ and
service times $B_p$. We thus get that
$\Psi_1(s)=-\Phi_A^{-1}(1/\Phi_{B_p}(s))$.
Since $\Phi_A(s)$ is strictly increasing in $s$, so is its inverse
$\Phi_A^{-1}(s)$. Combining this with  $\Phi_{B_p}(s)<\Phi_B(s)$
leads to the conclusion that $\Psi_1(s)<\Psi(s)$.
Recall that the residual busy period $P^*$ satisfies
$P^*\stackrel{d}{=} P(W)$; its decay rate is given by
$\sup_{s\in[0,\gamma_w]} \{s-\Psi(s)\}$, as can be seen by
mimicking the proof of Theorem 3.1.  Hence,
\beeq \label{pster} \gamma_{w_2}=\sup_{s\in [0,\gamma_w]} \{s-\Psi_1(s)\}>
\sup_{s\in[0,\gamma_w]} \{s-\Psi(s)\}=\gamma_{p^*}.\eneq
The proof is completed by recalling that
 for light-tailed service times, $\gamma_p= \gamma_{p^*}$.
\halmos\\

 From Theorem \ref{massinendpoint} and Proposition \ref{diffdecayrates} we conclude a similar result for the SRPT discipline.
\begin{Cor}\label{inbetween} If $0<\pr{B=x_B}<1,$ then $ \gamma_{p}< \gamma_{v}< \gamma_{w}.$
\end{Cor}
Hence, if there is mass in the endpoint $x_B$, then the decay rate of the sojourn time under SRPT lies strictly between those under LIFO and FIFO.\\

The following consequence of Theorem \ref{massinendpoint}  indicates that in some sense, all values between those of LIFO and FIFO are assumed.
Let  $F_q$ be the mixture  of a distribution bounded by $c$, and a distribution with all mass in $c$, such that  $\pr{B=c}=q$. Assume that $c<\expect A$, and let $\gamma_v(q)$ denote the decay rate
of the sojourn time in a queue with service-time distribution $F_q$.
\begin{Prop}  The decay rate  $\gamma_v(q)$ is continuous in $q$. In particular, it increases from
 $\gamma_p(0)$ to $\gamma_w(1)$, and assumes all values in between.
\end{Prop}
{\bf Proof}\\
By Theorem \ref{massinendpoint}, it is enough to show that $\gamma_w(q)$ and $\Psi_1(s)$ are continuous in $q$.
Since $\Phi_{B_1}$ is constant in $q$, and $\Phi_{A_1}$ is continuous in $q$, also $\Psi_1(s)$ is continuous in $q$.
Furthermore, since $\Phi_A$ is constant in $q$, $\Phi_B$ is continuous in $q$, and $\gamma_w(q)$ is finite for all $q$,  Proposition \ref{busyperiodprop}  implies that $\gamma_w(p)$ is continuous in $q$, and the proof is completed.\halmos\\

Table \ref{tabel} below shows the decay rates of the sojourn time that are known in the literature.
It turns out that the property in Corollary \ref{inbetween} is quite special, since almost all other known decay rates are either minimal or maximal.

\begin{table}[!h]
\begin{center}
 \begin{tabular}{|l|l|l|l|}
\hline
decay rate & discipline & condition & queue\\
\hline
\hline
 $\gamma_p$ (minimal), (\ref{eerste}) &  LCFS \cite{PR04}&  &$GI/GI/1$\\
&FB \cite{MN05}  &  & $M/GI/1$\\
& PS \cite{MZ05} &  $\forall c>0:  \log \pr{B>c\log x}=o(x)$ & $GI/GI/1$\\
&  ROS \cite{MZ05} &  &  $M/M/1$  \\
&  SRPT here & $\pr{B=x_B}=0$ & $GI/GI/1$\\
\hline
$\gamma$ (in between)  &  PS \cite{EZB05} & & $M/D/1$\\
$\gamma_{v}$  (in between), (\ref{drvsrpt}) & SRPT here & $0<\pr{B=x_B}<1$ & $GI/GI/1$\\
\hline
 $\gamma_w$ (maximal), (\ref{gammap}) & FCFS \cite{RS01}  & & $GI/GI/1$\\
& SRPT & $\pr{B=x_B}=1$ &  $GI/D/1$\\
\hline
\end{tabular}\\
\caption{\em The decay rate of the sojourn time under  several disciplines for light-tailed service times.}
\label{tabel}
\end{center}
\end{table}

\subsection{Poisson arrivals}\label{possonian}

Consider the priority queue of Section 3, with the additional
assumption that $A$ has an exponential distribution with rate $\lambda$.
Letting $\lambda_1=p\lambda$ and $\lambda_2=(1-p)\lambda$,
we  get that $\Psi_1(s)=\lambda_1 (\Phi_{B_1}(s)-1)$.
Thus, we have
\beeq
\gamma_{w_2}= \sup_{s\in [0,\gamma_w]} \{s-\lambda_1(\Phi_{B_1}(s)-1)\}.
\eneq
Suppose that $1- \lambda_1 \Phi_{B_1}'(\gamma_w)>0$.
Then the maximum value is attained in $\gamma_w$ and we have
\beeq \label{onze}
\gamma_{w_2}= \gamma_w- \lambda_1(\Phi_{B_1}(\gamma_w)-1).
\eneq
This expression is rather explicit, as $\gamma_w$ is the positive solution of
the equation
\beeq\label{defgammaw}\gamma_w=
\lambda (\Phi_B(\gamma_w)-1).
\eneq

The goal of this subsection is to show that in the case of Poisson arrivals,  
our  expression (\ref{onze}) coincides
with the expression of $\gamma_{w_2}$ given by Abate \& Whitt \cite{AW97}.
Assuming that $\expect{B_1}=1$, it is shown in \cite{AW97}, p.\ 18,  that
$-\gamma_{w_2}$ is the solution of
$\hat{f}(s)=1/\rho$, with
\[ \hat{f}(s)=\frac{\rho_1}{\rho_1+\rho_2} \hat{h}_0^{(1)}(s)+
\frac{\rho_2}{\rho_1+\rho_2}\hat{g}_{2e}(z_1(s)),
\qquad z_1(s)=s+\lambda_1-\lambda_1\hat{b}_1(s),\]
\[ \hat{h}_0^{(1)}(s)=\frac{1-\hat{b}_1(s)}{s+\rho_1-\rho_1\hat{b}_1(s)}=
\frac{1-\hat{b}_1(s)}{z_1(s)},
\qquad \hat{g}_{2e}(s)=\frac{1-\hat{g}_2(s)}{s g_{21}},\]
where $\hat{b}_1(s)$ is the LST of the $M/G/1$ busy period,
$\hat{g}_2(s)=\Phi_{B_2}(-s)$ and
$g_{21}=\expect{B_2}$.

Our expression of $\gamma_{w_2}$ seems preferable,
although we hasten to add that the form provided by \cite{AW97}
is more convenient when considering the more complicated task of obtaining
precise asymptotics, as is done in \cite{AW97}.

We now simplify the description of $\gamma_{w_2}$ in \cite{AW97}.
Since $\rho=\rho_1+\rho_2,$ we have
\[ \frac{1}{\rho_1+\rho_2}=\frac{1}{\rho}=\hat{f}(s)=
\frac{\rho_1}{\rho_1+\rho_2} \frac{1-\hat{b}_1(s)}{z_1(s)}+
\frac{\rho_2}{\rho_1+\rho_2}\frac{1-\hat{g}_2(z_1(s))}{z_1(s) \expect{ B_2}}.\]
Hence,
\[
z_1(s)=\rho_1[1-\hat{b}_1(s)]+ \frac{\rho_2}{\expect{B_2}}
[1-\hat{g}_2(z_1(s))].
\]
Consequently,
\[
s=\lambda_2 [1-\hat{g}_2(z_1(s))] = \lambda_2(1-\Phi_{B_2}(-z_1(s)).
\]
Since
$\lambda \Phi_B(s)=\lambda_1\Phi_{B_1}(s)+\lambda_2\Phi_{B_2}(s)$, we
can rewrite this into
\beeq \label{rewr1}
s +\lambda_1[1-\Phi_{B_1}(-z_1(s))] = \lambda [1-\Phi_{B}(-z_1(s))].
\eneq
Since the LST of the busy period satisfies the
fixed point equation
\beeq \label{fixedpoint}\hat{b}_1(s)= \Phi_{B_1}(-z_1(s)),\eneq
we can rewrite (\ref{rewr1}) as
\[  s+\lambda_1(1-\hat{b}_1(s))=\lambda [1-\Phi_{B}(-z_1(s))],\]
and thus, using the definition of $z_1(s)$,
\[
z_1(s)= \lambda [1-\Phi_{B}(-z_1(s))].
\]
Using the definition of $\gamma_w$ in (\ref{defgammaw}), we see that
$\gamma_{w_2}$ is the solution of
$\gamma_w=-z_1(s).$
We now give an alternative expression for $-z_1(s)$.
An alternative expression for the busy period transform was
found by Rosenkrantz \cite{Rosenkrantz83}:
defining $\phi(s)=\lambda_1(1-\Phi_{B_1}(-s))-s$,
it holds that
\beeq \label{altexp}
\hat{b}_1(s)= \Phi_{B_1}(-\phi^{-1}(s)).
\eneq
Since $\Phi_{B_1}(s)$ is strictly increasing, it follows from
(\ref{fixedpoint}) and (\ref{altexp}) that
$z_1(s)=\phi^{-1}(s)$. Hence, $\gamma_{w_2}$ is the solution of
$\gamma_w=-\phi^{-1}(s),$
and thus we obtain
\[
\gamma_{w_2}= \phi(-\gamma_w)= \gamma_w-\lambda_1(\Phi_{B_1}(\gamma_w)-1),
\]
which is indeed equal to our expression (\ref{onze}).\\

To conclude this section, we remark that for the M/G/1 queue, the decay rate $\gamma_v$ can take on  a simple form. Suppose that $1- \lambda_1 \Phi_{B_1}'(\gamma_w)>0$, and that $0<\pr{B=x_B}<1$.
Then by Theorem \ref{massinendpoint} and the expression for $\gamma_w$ given in (\ref{defgammaw}), we have
\begin{align*} \gamma_v& =\gamma_w- \lambda_1(\Phi_{B_1}(\gamma_w)-1)= \lambda (\Phi_B(\gamma_w)-1)- \lambda_1(\Phi_{B_1}(\gamma_w)-1)\\ & =\lambda_2(\Phi_{B_2}(\gamma_w)-1)= \lambda\pr{B=x_B} (e^{x_B\gamma_w}-1).
 \end{align*}

\subsection{Heavy traffic}\label{heavytraffic}

We now examine the behavior of the decay rate $\gamma_v$ of the SRPT sojourn time in heavy traffic. The aim of this section is to
show that the behavior of this decay rate critically depends upon whether $\pr{B=x_B}>0$ or not.
If $\pr{B=x_B}=0$, then $\gamma_v=\gamma_p$ by Theorem \ref{nomassinendpoint}. The results in Section 4.2 of \cite{MZ05} then imply
that $\gamma_v \sim C (1-\rho)^2$ for some constant $C$. We now show that a fundamentally different behavior
applies if $\pr{B=x_B}>0$.

Since, in this case, we have a relationship with the $GI/GI/1$ priority queue, we
consider first the setting of Section 3.  We let
 the service time $B_2$ increase in such a way that $\rho\rightarrow 1$.
Specifically, we consider a sequence of systems indexed by $r$,
such that $p$, $A_1$, $A_2$ and $B_1$ are all fixed, and that $B_2=B_2(r)$
is such that the traffic load satisfies $\rho_r=1-1/r$. Let $\gamma_w(\rho_r)$ denote
the decay rate of the workload in such a queue.

If we let $\sigma_A^2<\infty$ be the variance of $A$ and assume that the variance
 of $B(r)=B_1+B_2(r)$ converges to  $\sigma_B^2$, then
it holds that for $r\to\infty$ (cf.\ Corollary 3 of \cite{GW94}),
\begin{equation}\label{wht}
\gamma_w(\rho_r) \sim K(1-\rho_r),
\end{equation}
with $K=2/(\sigma_A^2+\sigma_B^2)$.
In particular, $\gamma_w(\rho_r) \downarrow 0$.
Consequently, if $r$ is large enough, we always have
$\gamma_{w_2}(\rho_r)= \gamma_{w}(\rho_r) -\Psi_1(\gamma_{w}(\rho_r))$ by Theorem \ref{w2thm}.
Since $\Psi_1(s) \sim \rho(1) s$ as $s\downarrow 0$, where $\rho(1)$ is the load in the high priority queue,
we obtain the following heavy-traffic result for $\gamma_{w_2}$.

\begin{Prop}
For $\rho\rightarrow 1$ as described above, we have
\[
\gamma_{w_2} \sim K(1-\rho(1)) (1-\rho).
\]
\end{Prop}

Thus, also $\gamma_v$ is of the order $(1-\rho)$ if $\pr{B=x_B}>0$.
This behavior is notably different from the  $(1-\rho)^2$ behavior of $\gamma_p$.

\section{Conditional sojourn times} \label{s6}

Our results in Section 4 and 5 show that the decay rate $\gamma_v$ for SRPT
is smaller than $\gamma_w$, which is the decay rate of the waiting
(and sojourn) time under FIFO.
Thus, one could say that according to this performance measure, SRPT
is worse than  FIFO.

The reason that the sojourn-time decay rate under SRPT is small is apparent when taking a
closer look at the proof in Section \ref{4.1}:
the sojourn time of a customer with a (very) large  service time  looks like a residual busy period.
However, smaller customers may have a much shorter sojourn time. In fact, for the {\em conditional}
sojourn time $V_{SRPT}(y) =[V_{SRPT} \mid B=y]$ under the preemptive SRPT discipline, the following proposition holds.

\begin{Prop} If  $\pr{B=y}=0$, then
$\log\pr{V_{SRPT}(y)>x} \sim -\gamma_p^y x$ as $x\to\infty$.
\end{Prop}
{\bf Proof}\\
For the lower bound, we remark that
$V_{SRPT}(y)$ is stochastically larger than the residual busy period $P^{*y}$ in the queue with service time $B^y$. This residual busy period has decay rate $\gamma_p^y$.
For the upper bound, we consider an alternative  queue with generic service time $B^y$, stationary workload at arrival instants  $W^y$ and busy period $P^y$. Now observe  that in the original queue, at any point in time, at most one customer with original service time larger than $y$ has remaining service time smaller han $y$. Hence, we can bound
\[ V_{SRPT}(y)\leq_{st} P^y(W^y+y+y),\]
where $P^y(x)$ is a busy period in the alternative queue starting with an exceptional customer of length $x$. Applying the Chernoff bound, and arguing like in the proof of Proposition \ref{busyperiodprop}, we find
\[ \limsup_{t\to\infty} \frac{1}{t} \log \pr{V_{SRPT}(y)>t}\leq - \sup_{s\in [0,\gamma_{w^y}]}\{s-\Psi^y(s)\}=-\gamma_{p*}^y,\]
where the last equality follows from (\ref{pster}).
The upper bound follows from noting that $P^y$ and $P^{*y}$ have the same decay rate, and  the proof is completed. \halmos\\

Suppose that $B$ has a density, so that $\gamma_p^y$ is continuous in $y$.
Then the  function $\gamma_p^y$ strictly decreases in  $y$, and converges to $\gamma_p>\gamma_w$ as $y\to\infty$. Further,  $\gamma_p^y\rightarrow\infty$ as $y\rightarrow 0$, since
$\Psi^y(s)\rightarrow 0$ as $y\rightarrow 0$.
Hence, there exists a critical value $y^*$ for which $\gamma_{p}^{y^*}=\gamma_w$.
Thus, when the decay rate is used as a performance measure, one
could say that FIFO is a better discipline than SRPT for customers
of size larger than $y^*$; the fraction of customers
that suffer from a change from FIFO to SRPT is $\pr{B>y^*}$.
We now describe the behavior of $y^*$ as a function of $\rho$ for $\rho\to 1$ and $\rho\to 0$.

\begin{Prop}  Let $y^*=\sup\{y: \gamma_p^y\geq \gamma_{w} \}.$ If $\rho \to 1$, then $y^*\to x_B$.
\end{Prop}
{\bf Proof}\\
 Let $y<x_B$ be fixed, and let $\gamma_p^y(\rho)$ be the decay rate of $P^y$ as a function of $\rho$, and define $\gamma_w(\rho)$ similarly.
Since $P^y$ is a busy period in a stable queue, even when $\rho=1$ in the original queue, we have  $\gamma_p^y(\rho)\geq \gamma_p^y(1)>0$ for all $\rho<1$.
By (\ref{wht}), we have for $\rho$ large enough,
\[ \gamma_w(\rho)<\gamma_p^y(1) \leq \gamma_p^y(\rho).\]
Hence, for $\rho$ large enough, $y^*\geq y$. Since $y<x_B$ was arbitrary, the proof is completed.\halmos

\begin{Prop}
If  the service time $B$ has decay rate $\gamma_b\in(0,\infty)$, then $y^*\to \infty$
for $\rho\to 0$.
\end{Prop}
{\bf Proof}\\
Let $\rho\to 0$ by setting the generic inter-arrival time equal to $rA$ and letting $r\to \infty$.
Since  $\Phi_{rA}(x)\to 0$ for all $x<0$, we have  $\Phi_{rA}^{-1}(x)\to 0$ for all $0<x<1$.
Hence, for all $y$,
\begin{equation}\label{rhois0}\gamma_p^y(\rho)=\sup_{s\geq 0} \Big\{s+ \Phi_{rA}^{-1}\Big(\frac{1}{\Phi_{B^y}(s)}\Big)\Big\}\to \infty, \qquad
\rho\to 0.\end{equation}
The workload does not depend on the discipline as long as the discipline is work-conserving.
Further, conditioned on it being positive, the workload under FIFO is stochastically larger than a residual service time, which for light-tailed distributions has the same decay rate as $B$.
Hence,  we have  $\gamma_w(\rho)\leq \gamma_b<\infty$ for all $\rho$. It then follows from (\ref{rhois0}) that $\gamma_w^y(\rho)> \gamma_w(\rho)$  eventually  as $\rho\to 0$ for all $y$, and we can conclude that  $y^*\to \infty$  as $\rho\to 0$.  \halmos

\subsection{Numerical example}\label{numsec}
\begin{figure}[h!]
\begin{center}
\epsfxsize=9cm 
\psfrag{r}{$\rho$}
\psfrag{t}{$\mathbb{P}\{B>y^*(\rho)\}$}
\epsfbox{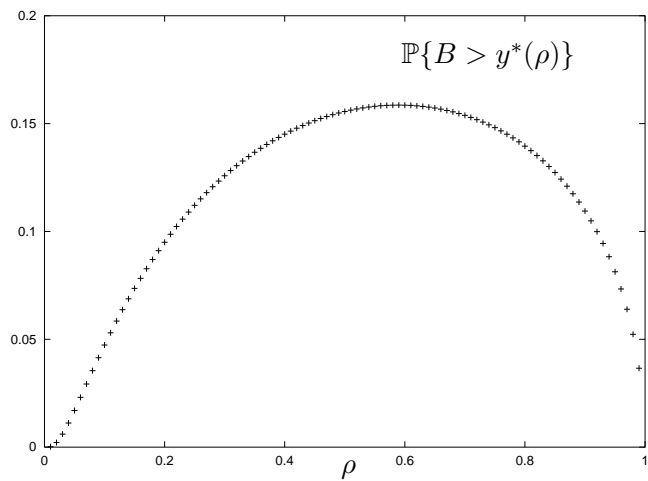}
\end{center}
\caption{The probabilities $\pr{B>y^*(\rho)}$ for $\rho\in (0,1)$ in the $M/M/1$ queue.}
\label{figyster}
\end{figure}

As an illustration, we compute $y^*$ and $\pr{B>y^*}$
for the $M/M/1$ queue with $\expect{B}=1$ and arrival rate $\lambda$ (so that $\rho=\lambda$). Figure \ref{figyster} shows the probabilities  $\pr{B>y^*}$
for various values of $\rho$.  \\

\noindent
 From the figure, it is clear that $y^*$ becomes very large under low and high loads. But even for moderate values of $\rho$ it is clear that
about 85 percent of the customers would prefer (from a large-deviations point of view) SRPT over FIFO.

\section{Conclusions}
\label{s7}
To conclude the paper, we summarize our results. For the $GI/GI/1$ queue with light-tailed service times, we obtained expressions for the logarithmic
decay rate of the tail of the workload, the busy period, the waiting time and sojourn time of low-priority customers in a priority queue, and the sojourn time under the
(preemptive and non-preemptive) SRPT discipline.

For the sojourn time under SRPT, it turns out that there are three different regimes, namely for service times with no mass, with some mass and  with all mass in the endpoint of the service-time distribution. In the first case the decay rate is minimal among all work-conserving disciplines, in the last case it is maximal, but if there is some mass in the endpoint, then the decay rate lies strictly in between these two.
The large-deviations results for the unconditional sojourn times suggest that a switch from FIFO to SRPT is not advisable.
The results in Section 6 show that this suggestion is only valid for very large service times: in the $M$/$M$/1 queue, at least about 85 percent of the customers would benefit from a change from FIFO to SRPT.

There are several topics that  are interesting for further research. First of all,  large deviations for the queue length under SRPT are not well understood. 
A second problem is to obtain precise asymptotics for the tail behavior of the low-priority waiting time, or perhaps even the sojourn time. Finally, it would be interesting to compare conditional sojourn times of  FIFO and PS from a large-deviations
point of view. It is not clear to us which discipline performs better, and what the influence of the job size might be.\\

\noindent
{\bf Acknowledgments}\\
We would like to thank Marko Boon
 for helping us out with the numerics in Section \ref{numsec}, and Ton Dieker and Michel Mandjes for several useful
comments.

\end{document}